\documentclass{article}[12pt]
\usepackage{amssymb}
\usepackage{epsfig}
\usepackage{amsmath}
\usepackage{amsfonts}
\usepackage{pgf,tikz}

\newtheorem{theorem}{\sc Theorem}[section]

\newtheorem{definition}[theorem]{\sc Definition}
\newtheorem{remark}[theorem]{\sc Remark}

\def\qed{\hbox to 0pt{}\hfill$\rlap{$\sqcap$}\sqcup$\medbreak}

\title{Fourier series for nonperiodic functions}

\author{Rodrigo L\'opez Pouso\footnote{Partially supported by
Agencia Estatal de Investigaci\'on, Spain, project PID2020-113275GB-I00.}}

\date{}
\begin{document}
 \maketitle

\begin{center}  {CITMAga, 15782, Santiago de Compostela, Spain\\ and \\ Departamento de An\'alise Matem\'atica, Estat\'{\i}stica e Optimizaci\'on \\Universidade de Santiago de Compostela \\ 15782, Faculty of Mathematics, Campus Vida\\
Santiago de Compostela, Spain.\\  Email: rodrigo.lopez@usc.es}
\end{center}

\medbreak

 \abstract{We introduce a small change in the definition of the Fourier series so that we can guarantee the coincidence with the given function at the endpoints of the interval even if the function does not assume the same value at the endpoints. This definition of the Fourier series also wipes out the Gibbs phenomenom at the endpoints of the interval and proves useful in the resolution of antiperiodic boundary value problems with linear partial differential equations.}

\medbreak


\medbreak

\section{Introduction}

The Fourier series of a function $f:[-\pi,\pi] \longrightarrow \mathbb R$ is
\begin{equation}
\label{sf}
Sf(x)=\dfrac{a_0}{2}+\sum_{n=1}^{\infty} \left( a_n \cos nx +b_n \sin nx
\right), \quad x \in \mathbb R,
\end{equation}
where
\begin{eqnarray}
\label{an}
a_n=\dfrac{1}{\pi}\int_{-\pi}^{\pi}{f(x) \cos nx \, dx}, \, \, n=0,1,2,\dots,\\
\label{bn}
b_n= \dfrac{1}{\pi}\int_{-\pi}^{\pi}{f(x) \sin nx \, dx}, \, \, n=1,2,3,\dots,
\end{eqnarray}
and provided that the integrals exist (they do exist if $f$ is just Lebesgue integrable on $[-\pi,\pi]$, although the series $Sf(x)$ need not converge for every $x \in \mathbb R$).

We refer readers to \cite{m, s} for the basic results on Fourier series and their applications.

\bigbreak

For instance, the Fourier series of $f(x)=x$, $x \in [-\pi,\pi]$, is
$$Sf(x)=2 \sum_{n=1}^{\infty}\dfrac{(-1)^{n+1}}{n} \sin nx, \, \, x \in \mathbb R.$$

This simple example shows two disadvantages of the Fourier series which we overcome here. First, although $f= Sf$ on $(-\pi,\pi)$, we have $f(\pm \pi) \neq Sf(\pm \pi)=0$ (see Figure 1); second, the partial sums of $Sf(x)$ exhibit the so--called Gibbs phenomenon around $\pm \pi$, i.e. the error of approximating $f(x)$ by partial sums of $Sf(x)$ is greater around $\pm\pi$ than around any other point in $(-\pi,\pi)$, see Figure 2.

\begin{figure}[h]
\centering
\includegraphics[width=7cm,keepaspectratio=true]{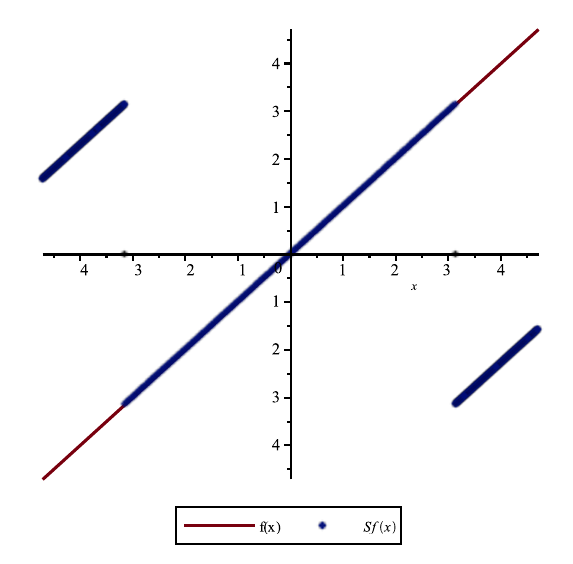}
\caption{Observe that $f(\pm \pi) \neq Sf(\pm \pi)=0$.}
\end{figure}

\begin{figure}[h]
\centering
\includegraphics[width=7cm,keepaspectratio=true]{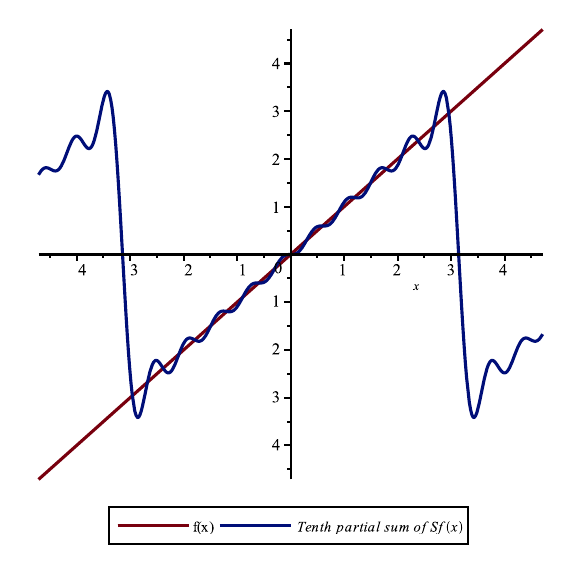}
\caption{Observe the Gibbs phenomenom around $\pm \pi$.}
\end{figure}

These two drawbacks essentially come from the fact that $f(x)=x$ is not periodic on $[-\pi,\pi]$, in the sense\footnote{We will say $f$ is periodic on $[-\pi,\pi]$ if $f(-\pi)=f(\pi)$, but this is an abuse of language because $f$ is not defined on the whole of $\mathbb R$, and therefore it cannot be periodic in the strict sense. However, notice that if $f(-\pi)=f(\pi)$, then $f$ can be extended as a $2\pi$--periodic function on $\mathbb R$. A similar abuse of language will occur in Section 2 when defining antiperiodic functions.} that $f(-\pi) \neq f(\pi)$, and Fourier series are always periodic.
\bigbreak

There exist many convergence criteria for Fourier series. The following result, whose proof can be looked up in \cite[Corollary 8.48]{s}, is especially appropriate for the purposes of this paper.

\begin{theorem}
\label{tu}
Let $f:[-\pi,\pi] \longrightarrow \mathbb R$ be such that  $f(-\pi)=f(\pi)$. 

If $f$ is continuous and has bounded variation on $[-\pi,\pi]$, then $Sf(x)$ is uniformly convergent on $[-\pi,\pi]$ and $Sf(x)=f(x)$ for all $x \in [-\pi,\pi]$.
\end{theorem}

Remember that continuity alone is not enough for the convergence of $Sf(x)$, see \cite[Remark 8.26]{s}.

On the other hand, recall that any of the following conditions implies that $f$ is continuous and has bounded variation on $[-\pi,\pi]$:
\begin{enumerate}
\item[(i)] $f$ is absolutely continuous on $[-\pi,\pi]$; 
\item[(ii)] $f$ Lipschitz continuous on $[-\pi,\pi]$;
\item[(iii)] $f$ is continuous and there is a partition
$-\pi=x_0 < x_1<\cdots<x_p=\pi,$ such that the restrictions $f_{|[x_{k-1},x_k]}$ are continuously differentiable on $[x_{k-1},x_k]$ for each $k=1,2,\dots,p$;
\item[(iv)] $f$ is continuously differentiable on $[-\pi,\pi]$.
\end{enumerate}

 Of course, our first example $f(x)=x$, $x \in [-\pi,\pi]$, is not periodic, so it falls outside the scope of Theorem 1.1.  

Can we have an alternative definition of $Sf$ so that Theorem \ref{tu} holds true if we remove the condition that $f(-\pi)=f(\pi)$? The answer is positive and the alternative definition of $Sf$ is pretty simple: we only have to replace the basic trigonometric functions $\cos nx$ and $\sin nx$ by $\cos (nx+x/2)$ and $\sin(nx+x/2)$.

\section{Fourier series for antiperiodic and nonperiodic functions}
We start with the somewhat opposite case to periodic functions, namely, that of antiperiodic functions. Then we will reduce the general (nonperiodic) case to the antiperiodic case. 

Let $f:[-\pi,\pi] \longrightarrow \mathbb R$ be continuous and have bounded variation on $[-\pi,\pi]$. Assume also that $f$ is antiperiodic in the mild sense that $f(-\pi)=-f(\pi)$.

The fundamental ideas in this paper lean on the functions 
$$f_c(x)=f(x) \cos \dfrac{x}{2} \quad \mbox{and} \quad f_s(x)=f(x)\sin \dfrac{x}{2}, \quad x \in [-\pi,\pi],$$ 
which are continuous and have bounded variation on $[-\pi,\pi]$ and, moreover, they are periodic because $f$ is antiperiodic, i.e. $f_c(-\pi)=f_c(\pi)$ and $f_s(-\pi)=f_s(\pi)$. 

Hence, Theorem \ref{tu} can be applied to $f_c$ and $f_s$ and we get
\begin{equation}
\label{sfc}
f_c(x)=Sf_c(x)=\dfrac{a_0}{2}+\sum_{n=1}^{\infty} \left( a_n \cos nx +b_n \sin nx
\right), \quad x \in [-\pi,\pi],
\end{equation}
\begin{equation}
\label{sfs}
f_s(x)=Sf_s(x)=\dfrac{\tilde{a}_0}{2}+\sum_{n=1}^{\infty} \left( \tilde{a}_n \cos nx +\tilde{b}_n \sin nx
\right), \quad x \in [-\pi,\pi],
\end{equation}
where, according to (\ref{an}) and (\ref{bn}),
\begin{eqnarray}
\label{anc}
a_n=\dfrac{1}{\pi}\int_{-\pi}^{\pi}{f_c(x) \cos nx \, dx}=\dfrac{1}{\pi}\int_{-\pi}^{\pi}{f(x) \cos \dfrac{x}{2} \cos nx \, dx}, \, \, n=0,1,2,\dots,\\
\label{bnc}
b_n= \dfrac{1}{\pi}\int_{-\pi}^{\pi}{f_c(x) \sin nx \, dx}=\dfrac{1}{\pi}\int_{-\pi}^{\pi}{f(x) \cos \dfrac{x}{2} \sin nx \, dx}, \, \, n=1,2,3,\dots,
\end{eqnarray}
and \begin{eqnarray}
\label{ans}
\tilde{a}_n=\dfrac{1}{\pi}\int_{-\pi}^{\pi}{f_s(x) \cos nx \, dx}=\dfrac{1}{\pi}\int_{-\pi}^{\pi}{f(x) \sin \dfrac{x}{2} \cos nx \, dx}, \, \, n=0,1,2,\dots,\\
\label{bns}
\tilde{b}_n= \dfrac{1}{\pi}\int_{-\pi}^{\pi}{f_s(x) \sin nx \, dx}=\dfrac{1}{\pi}\int_{-\pi}^{\pi}{f(x) \sin \dfrac{x}{2} \sin nx \, dx}, \, \, n=1,2,3,\dots
\end{eqnarray}

Obviously, for all $x \in [-\pi,\pi]$ we have
 $$f(x)=f_c(x) \cos \dfrac{x}{2}+f_s(x) \sin \dfrac{x}{2},$$
 and we use (\ref{sfc}) and (\ref{sfs}) to obtain that
\begin{align}
\nonumber
f(x)&=\dfrac{a_0}{2} \cos \dfrac{x}{2}+\dfrac{\tilde{a}_0}{2} \sin \dfrac{x}{2}  
 +\sum_{n=1}^{\infty}\left(a_n \cos nx \cos \dfrac{x}{2}+b_n \sin nx \cos \dfrac{x}{2} \right)   \\
 \nonumber
 &\quad +\sum_{n=1}^{\infty}\left(\tilde {a}_n \cos nx\sin \dfrac{x}{2}+\tilde {b}_n \sin nx \sin \dfrac{x}{2} \right)  \\
 \nonumber
 &=\dfrac{a_0}{2} \cos \dfrac{x}{2}+\dfrac{\tilde{a}_0}{2} \sin \dfrac{x}{2}  \\
 \label{suma}
 &\quad +\lim_{m \to \infty} \sum_{n=1}^{m}\left( 
 a_n \cos nx \cos \dfrac{x}{2}+\tilde {a}_n \cos nx  \sin \dfrac{x}{2} 
+ b_n \sin nx\cos \dfrac{x}{2}+\tilde {b}_n \sin nx \sin \dfrac{x}{2}
 \right).
 \end{align}

Now we use the identities
\begin{eqnarray*}
2\cos nx \, \cos (x/2)=\cos \left(nx+x/2 \right)+\cos \left(nx-x/2 \right),\\
2\cos nx \, \sin  (x/2)=\sin \left(nx+x/2 \right)-\sin \left(nx-x/2 \right),\\
2\sin nx \, \cos (x/2)=\sin \left(nx+x/2 \right)+\sin \left(nx-x/2 \right),\\
2  \sin  nx \, \sin  (x/2)=\cos \left(nx-x/2 \right)-\cos \left(nx+x/2 \right),
\end{eqnarray*}
and we obtain the following expression for each fixed $m \in \mathbb N$:
\begin{align*}
\sum_{n=1}^{m}&\left(  a_n \cos nx \cos \dfrac{x}{2}+\tilde {a}_n \cos nx  \sin \dfrac{x}{2} 
+ b_n \sin nx\cos \dfrac{x}{2}+\tilde {b}_n \sin nx \sin \dfrac{x}{2}
 \right)\\
&=\sum_{n=1}^{m} \dfrac{a_n-\tilde{b}_n}{2} \cos \left( nx+\dfrac{x}{2} \right)+\sum_{n=1}^{m}\dfrac{\tilde{a}_n+b_n}{2} \sin \left( nx+\dfrac{x}{2} \right)\\
& \quad +\sum_{n=1}^{m} \dfrac{a_n+\tilde{b}_n}{2} \cos \left( nx-\dfrac{x}{2} \right)+\sum_{n=1}^{m}\dfrac{b_n-\tilde{a}_n}{2} \sin \left( nx-\dfrac{x}{2} \right)\\
&=\sum_{n=1}^{m} \dfrac{a_n-\tilde{b}_n}{2} \cos \left( nx+\dfrac{x}{2} \right)+\sum_{n=1}^{m}\dfrac{\tilde{a}_n+b_n}{2} \sin \left( nx+\dfrac{x}{2} \right)\\
& \quad +\dfrac{a_1+\tilde{b}_1}{2}\cos \dfrac{x}{2}+\dfrac{b_1-\tilde{a}_1}{2}\sin \dfrac{x}{2}\\
& \quad +\sum_{n=1}^{m} \dfrac{a_{n+1}+\tilde{b}_{n+1}}{2} \cos \left( nx+\dfrac{x}{2} \right)+\sum_{n=1}^{m}\dfrac{b_{n+1}-\tilde{a}_{n+1}}{2} \sin \left( nx+\dfrac{x}{2} \right).
\end{align*}
We return to (\ref{suma}) and we deduce that for all $x \in [-\pi,\pi]$ we have
\begin{align*}
f(x)&=\dfrac{a_0+a_1+\tilde{b}_1}{2}\cos \dfrac{x}{2}+\dfrac{\tilde{a}_0-\tilde{a}_1+b_1}{2} \sin \dfrac{x}{2}\\
& \qquad +\lim_{m \to \infty}  \sum_{n=1}^{m}\left(
  \dfrac{a_n +a_{n+1}-\tilde {b}_{n}+\tilde {b}_{n+1}}{2}   \cos \left( nx+\dfrac{x}{2} \right) +  \dfrac{\tilde{a}_n -\tilde{a}_{n+1}+{b}_{n}+{b}_{n+1}}{2}  \sin \left(nx+\dfrac{x}{2} \right) \right).
\end{align*}
Finally, we use the definitions (\ref{anc}), (\ref{bnc}), (\ref{ans}) and (\ref{bns}) and we arrive at the desired alternative Fourier series expansion for $f$: for all $x \in [-\pi, \pi]$ we have
\begin{equation}
\label{asf1}
f(x)=\sum_{n=0}^{\infty}\left( \alpha_n \cos \left(nx+\dfrac{x}{2} \right)  + \beta_n  \sin  \left(nx+\dfrac{x}{2} \right) \right),
\end{equation}
where
 $$\alpha_n=\dfrac{1}{\pi}\int_{-\pi}^{\pi}{f(x) \cos \left(nx+\dfrac{x}{2} \right) \, dx}, \, \, \beta_n=\dfrac{1}{\pi}\int_{-\pi}^{\pi}{f(x) \sin \left(nx+\dfrac{x}{2} \right) \, dx} \quad (n=0,1,2,3,\dots).$$
 Moreover, the series in (\ref{asf1}) converges uniformly on $[-\pi,\pi]$ thanks to uniform convergence of $Sf_c$ and $Sf_s$ on $[-\pi,\pi]$.
 
 Let us define the antiperiodic Fourier series of $f$ as
\begin{equation*}
AS f(x)=\sum_{n=0}^{\infty}\left( \alpha_n \cos \left(nx+\dfrac{x}{2} \right)+ \beta_n \sin  \left(nx+\dfrac{x}{2} \right) \right).
\end{equation*}

We have proven the following result, which is the analogue to Theorem \ref{tu} for anti--periodic functions. 

\begin{theorem}
\label{t1}
Let $f:[-\pi,\pi] \longrightarrow \mathbb R$ be an antiperiodic function, i.e. $f(-\pi)=-f(\pi)$. 

If $f$ is continuous and has bounded variation on $[-\pi,\pi]$, then the antiperiodic Fourier series $AS f(x)$ is uniformly convergent on $[-\pi,\pi]$ and $AS f(x)=f(x)$ for all $x \in [-\pi,\pi]$.

\end{theorem}

For $f(x)=x$, $x \in [-\pi,\pi]$, the antiperiodic Fourier series is
$$AS\, f(x)=\dfrac{8}{\pi}\sum_{n=0}^{\infty}\dfrac{(-1)^n}{(2n+1)^2}\sin  \left(nx+\dfrac{x}{2} \right),$$
which agrees with $f$ also at $\pm \pi$ and does not exhibit the Gibbs phenomenon around these points. See Figures 3 and 4.
\begin{figure}[h]
\centering
\includegraphics[width=7cm,keepaspectratio=true]{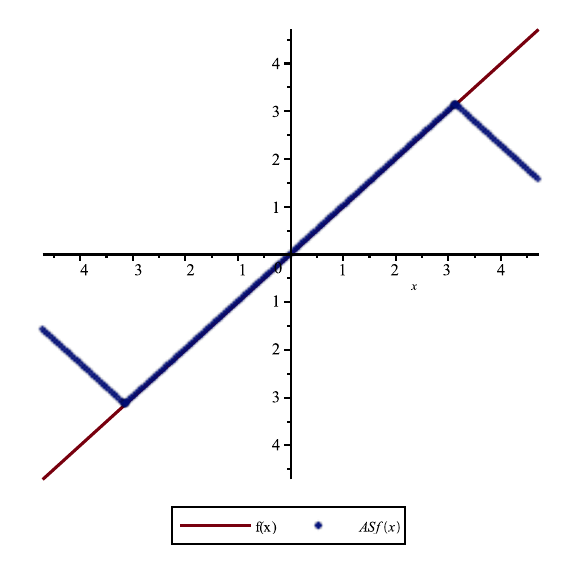}
\caption{Observe that $f(\pm \pi) =ASf(\pm \pi)$.}
\end{figure}

\begin{figure}[h]
\centering
\includegraphics[width=7cm,keepaspectratio=true]{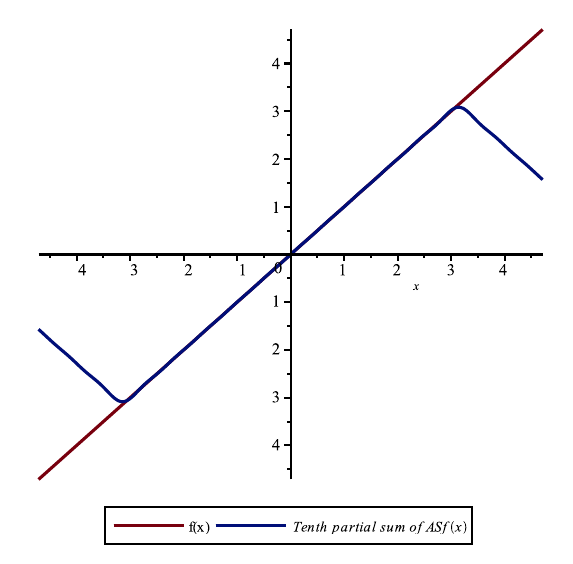}
\caption{No Gibbs phenomenom around $\pm \pi$.}
\end{figure}

Let us consider now the general nonperiodic case.

The main idea is that any function can be shifted upwards or downwards to become antiperiodic, so the antiperiodic Fourier series is useful regardless of the behaviour of $f$ at the boundary.

Indeed, let $f:[-\pi,\pi] \longrightarrow \mathbb R$ be an arbitrary function. It suffices to take $\gamma \in \mathbb R$ so that $F(x)=f(x)-\gamma$ be antiperiodic, i.e. we define
$$F(x)=f(x)-\gamma \quad (x \in [-\pi,\pi]) \quad \mbox{for} \quad \gamma=\dfrac{f(-\pi)+f(\pi)}{2},$$
and $F$ is antiperiodic.

Now, if $f$ is continuous and has bounded variation on $[-\pi,\pi]$ then so $F$ is, and then Theorem \ref{t1} ensures that 
$$f(x)-\gamma=F(x)=AS F(x)=\sum_{n=0}^{\infty}\left( \alpha_n \cos \left(nx+\dfrac{x}{2} \right)  + \beta_n  \sin  \left(nx+\dfrac{x}{2} \right)\right), \quad x \in [-\pi,\pi],$$
where
 $$\alpha_n=\dfrac{1}{\pi}\int_{-\pi}^{\pi}{(f(x)-\gamma) \cos \left(nx+\dfrac{x}{2} \right) \, dx}, \, \, \beta_n=\dfrac{1}{\pi}\int_{-\pi}^{\pi}{(f(x)-\gamma) \sin \left(nx+\dfrac{x}{2} \right) \, dx} \quad (n=0,1,2,3,\dots).$$
 Hence,
 $$f(x)=\gamma+\sum_{n=0}^{\infty}\left( \alpha_n \cos \left(nx+\dfrac{x}{2} \right)  + \beta_n  \sin  \left(nx+\dfrac{x}{2} \right)\right), \quad x \in [-\pi,\pi].$$
 
Finally, we put everything in a definition and a theorem. For the reader's convenience, we rescale formulas so that our results apply to intervals of the form $[-L,L]$, $L>0$.

\begin{definition}
\label{def}
The antiperiodic Fourier series of a function $f:[-L,L] \longrightarrow \mathbb R$, $L>0$, is
$$AS f(x)=\gamma+\sum_{n=0}^{\infty}\left( \alpha_n \cos \left(\dfrac{n \pi x}{L}+\dfrac{\pi x}{2L} \right)  + \beta_n  \sin  \left(\dfrac{n \pi x}{L}+\dfrac{\pi x}{2L} \right)\right),$$
where $\gamma=(f(-L)+f(L))/2$ and for $n=0,1,2,3,\dots$
 $$\alpha_n=\dfrac{1}{L}\int_{-L}^{L}{(f(x)-\gamma) \cos \left(\dfrac{n \pi x}{L}+\dfrac{\pi x}{2L} \right) \, dx}, \, \, \beta_n=\dfrac{1}{L}\int_{-L}^{L}{(f(x)-\gamma) \sin \left(\dfrac{n \pi x}{L}+\dfrac{\pi x}{2L} \right) \, dx},$$
provided that the integrals exist.
\end{definition}
 \begin{theorem}
\label{t1}
Let $f:[-L,L] \longrightarrow \mathbb R$ be continuous and have bounded variation on $[-L,L]$, $L>0$. 
 
Then the antiperiodic Fourier series $AS f(x)$ introduced in Definition \ref{def} is uniformly convergent on $[-L,L]$ and $AS f(x)=f(x)$ for all $x \in [-L,L]$.
\end{theorem}

As an example, we consider $f(x)=x^2+2x+1$ for $x \in [-1,1]$. Its antiperiodic Fourier series is 
$$AS f(x)=2-\dfrac{32}{\pi^3}\sum_{n=0}^{\infty}\dfrac{(-1)^n}{(2n+1)^3} \cos \left(n\pi x+\dfrac{\pi x}{2} \right)  +\dfrac{16}{\pi^2} \sum_{n=0}^{\infty}\dfrac{(-1)^n}{(2n+1)^2}  \sin  \left(n\pi x+\dfrac{\pi x}{2} \right),$$
and, according to Theorem \ref{t1}, we have $f(x)=AS f(x)$ for all $x \in [-1,1]$. The usual Fourier series is
$$Sf(x)=\dfrac{4}{3}+\dfrac{4}{\pi^2}\sum_{n=0}^{\infty}\dfrac{(-1)^n}{n^2} \cos \left(n\pi x\right)  -\dfrac{4}{\pi} \sum_{n=0}^{\infty}\dfrac{(-1)^n}{n}  \sin  \left(n\pi x \right).$$
Observe that, unlike the usual Fourier series, the antiperiodic Fourier series does not exhibit the Gibbs phenomenom at $\pm 1$, see Figures 5 and 6. Moreover, it is clear from the expressions of $AS f$ and $Sf$, or from Figures 5 and 6, that the antiperiodic Fourier series has a faster rate of convergence than the usual Fourier series. This is not a surprise, because the antiperiodic Fourier series uses regular Fourier series for periodic functions (namely, the functions $f_c$ and $f_s$ that we introduced at the beginning of this section) which accelerates convergence.

\begin{figure}[h]
\centering
\includegraphics[width=7cm,keepaspectratio=true]{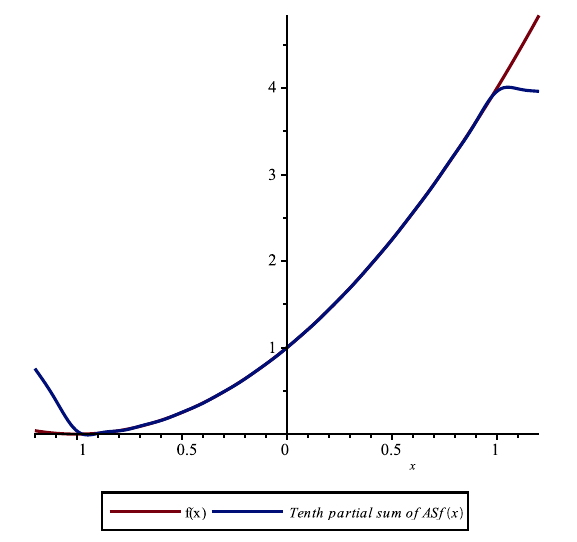}
\caption{No Gibbs phenomenom around $\pm 1$ and higher rate of convergence.}
\end{figure}

\begin{figure}[h]
\centering
\includegraphics[width=7cm,keepaspectratio=true]{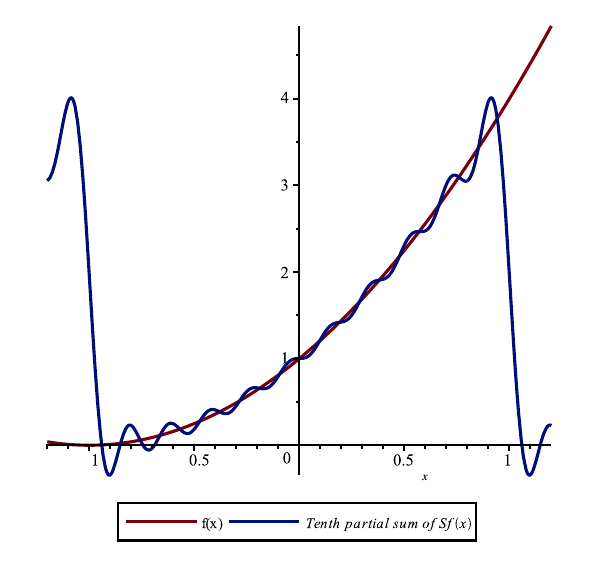}
\caption{Gibbs phenomenom around $\pm 1$ and lower rate of convergence.}
\end{figure}

We have only adapted Theorem \ref{tu} to the new setting of antiperiodic Fourier series, but surely every result on classical Fourier series has its corresponding analogue for antiperiodic Fourier series. To illustrate this possibility, and due to the importance of the theorem itself, we show how to prove Jordan's convergence test for the antiperiodic Fourier series as a consequence of the same result for the usual Fourier series. Further theorems on the Fourier series can be adapted for the antiperiodic Fourier series through analogous arguments.

\begin{theorem}
(Jordan's convergence test) Let $f \in L^1(-L,L)$, $L>0$, and let $\bar f$ denote its $2L$--periodic extension to $\mathbb R$, i.e. $\bar f=f$ a.e. on $(-L,L)$ and $\bar f(x+2L)=\bar f(x)$ for a.e. $x \in \mathbb R$.

Let $x_0 \in \mathbb R$ be fixed. If $\bar f$ is well defined and it has bounded variation on a neighborhood of $x_0$, then the antiperiodic Fourier series of $f$ converges for $x=x_0$ and 
$$AS f(x_0)=\dfrac{\bar f (x_0^+)+\bar f (x_0^-)}{2}.$$
\end{theorem}

\noindent
{\bf Proof.} We assume, without loss of generality, that $f$ is antiperiodic (so $\gamma=0$ in Definition \ref{def}). Adjusting the ideas by means of which we deduced the identity (\ref{asf1}), we can prove that
\begin{equation}
\label{asf2}
AS f(x)=Sf_c(x) \cos \dfrac{\pi x}{2L}+Sf_s(x) \sin \dfrac{\pi x}{2L}, \quad x \in \mathbb R,
\end{equation}
where $Sf_c$ and $Sf_s$ are the usual Fourier series of 
$$f_c(x)=f(x) \cos \dfrac{\pi x}{2L} \quad \mbox{and} \quad f_s(x)=f(x) \sin \dfrac{\pi x}{2L}, \quad \mbox{for a.a. $x \in [-L,L]$.}$$
The validity of (\ref{asf2}) depends on the convergence of $Sf_c(x )$ and $Sf_s(x )$, which we prove below for $x=x_0$.

Let $\bar f_c$ and $\bar f_s$ denote the $2L$--periodic extensions of $f_c$ and $f_s$, respectively. The assumptions imply that $\bar f_c$ and $\bar f_s$ satisfy the conditions in Jordan's test (see \cite[Theorem 8.47]{s}) and therefore $Sf_c(x_0)$ and $Sf_s(x_0)$ are convergent and
\begin{align*}
AS f(x_0)&=Sf_c(x_0) \cos \dfrac{\pi x_0}{2L}+Sf_s(x_0) \sin \dfrac{\pi x_0}{2L} \\
&=\dfrac{\bar f_c(x_0^+)+\bar f_c (x_0^-)}{2}\cos \dfrac{\pi x_0}{2L}+\dfrac{\bar f_s(x_0^+)+\bar f_s (x_0^-)}{2}\sin \dfrac{\pi x_0}{2L}\\
&=\dfrac{\bar f (x_0^+)+\bar f (x_0^-)}{2}.
\end{align*}
\qed

\begin{remark}
The following sufficient condition is very useful in applications of the Jordan's test: if $\bar f$ is differentiable on $(x_0-r,x_0) \cup (x_0,x_0+r)$ for some $r>0$ and there exist
$$\lim_{x \to x_0^{\pm}} \bar f'(x) \in \mathbb R \quad \mbox{(they need not be equal),}$$
then $\bar f$ has bounded variation on a neighborhood of $x_0$.
\end{remark}

\section{Application to PDEs and separation of variables}
Obviously, the antiperiodic Fourier series can be used to solve PDE problems with antiperiodic boundary conditions, but not only these ones.

As a sample, consider the following problem with the heat equation:
\begin{eqnarray}
\label{edp}
u_{t}-k \, u_{xx}=0, \quad (x,t) \in (-\pi,\pi)\times (0,\infty),\\
\label{ic}
u(x,0)=f(x), \quad x \in [-\pi,\pi],\\
\label{b1}
u(-\pi,t)+u(\pi,t)=2, \quad t  \ge 0,\\
\label{b2}
u_x(-\pi,0)+u_x(\pi,0)=0, \quad t \ge 0.
\end{eqnarray}
 
We assume that $k>0$ and $f:[-\pi,\pi] \longrightarrow \mathbb R$ is continuously differentiable on $[-\pi,\pi]$. We also assume that $f$ is compatible with the boundary conditions, i.e. $f(-\pi)+f(\pi)=2$ and there exist $f'(\pm \pi)$ and $f'(-\pi)+f'(\pi)=0$.

This problem is a mathematical model for the evolution of the temperature $u(x,t)$ at the point $x$ and time $t$ of a thin rod that we identify with the interval $[-\pi,\pi]$. For physical reasons, we should assume that $f(x) \ge 0$ for all $x \in [-\pi,\pi]$, because $f(x)$ is the temperature of the point $x$ at the initial time $t=0$. 

The boundary conditions can be interpreted as follows. Condition (\ref{b1}) means that the mean value of the temperature at the endpoints is kept constantly equal to $1$. In turn, condition (\ref{b2}) means that the heat flux at the endpoints is equal in absolute value but it has opposite sign, so we either have heat entering the rod from both ends at the same rate or leaving the rod from both ends at the same rate.

\bigbreak

One can prove, either by energy methods as in \cite{m} or by Fourier series as in \cite{p}, that (\ref{edp})--(\ref{b2}) has a unique solution $u \in {\cal C}^1([-\pi,\pi]\times [0,\infty))\cap {\cal C}^2((-\pi,\pi) \times (0,\infty))$. However, we are more interested in showing how the previous theory can be used in combination with the method of separation of variables for the computation of the solution. 

Observe that (\ref{b1}) is not an homogeneous boundary condition, so we cannot use the method directly. Clearly, the solution can be obtained as
\begin{equation}
\label{solu}
u(x,t)=v(x,t)+1, \quad (x,t) \in [-\pi,\pi] \times [0,\infty),
\end{equation}
where $v(x,t)$ is the solution of the following problem with homogeneous boundary data:
\begin{eqnarray}
\label{edp2}
v_{t}-k \, v_{xx}=0, \quad (x,t) \in (-\pi,\pi)\times (0,\infty),\\
\label{ic2}
v(x,0)=f(x)-1, \quad x \in [-\pi,\pi],\\
\label{b11}
v(-\pi,t)+v(\pi,t)=0, \quad t  \ge 0,\\
\label{b22}
v_x(-\pi,0)+v_x(\pi,0)=0, \quad t \ge 0.
\end{eqnarray}

Now we use separation of variables to compute $v(x,t)$. As usual, we start by computing nontrivial solutions of (\ref{edp2}), (\ref{b11}) and (\ref{b22}) which have the form $v(x,t)=X(x)T(t)$. Inserting this $v(x,t)$ in  (\ref{edp2}), (\ref{b11}) and (\ref{b22}) leads to finding all possible values of $\lambda \in \mathbb R$ for which we have nontrivial solutions of
\begin{equation}
\label{x}
X''(x)=\lambda X(x), \quad X(-\pi)=-X(\pi), \, \, X'(-\pi)=-X'(\pi),
\end{equation}
and
\begin{equation}
\label{t}
T'(t)=\lambda T(t), \, \, t \ge 0.
\end{equation}
It is easy to prove that (\ref{x}) has only the zero solution if $\lambda \ge 0$. For $\lambda <0$, the general solution of the ODE in (\ref{x}) is
\begin{equation}
\label{solgral}
X(x)=A \cos \left( \sqrt{-\lambda} x\right )+B \sin  \left( \sqrt{-\lambda} x\right ), \quad A, B \in \mathbb R,
\end{equation}
and then we impose the boundary conditions in (\ref{x}) to get
$$X(-\pi)=-X(\pi) \Leftrightarrow 2A \cos  \left( \sqrt{-\lambda} \pi\right )=0,$$
and
$$X'(-\pi)=-X'(\pi) \Leftrightarrow 2B \cos  \left( \sqrt{-\lambda} \pi\right )=0.$$
Therefore, if $\cos \left( \sqrt{-\lambda} \pi\right ) \neq 0$, then we have $A=B=0$, so we only get the zero solution. If, on the other hand, we have $\cos \left( \sqrt{-\lambda} \pi\right )=0$, then (\ref{solgral}) solves (\ref{x}) for any values of $A$ and $B$. Observe that $\sqrt{-\lambda} \pi >0$, so 
$$\cos \left( \sqrt{-\lambda} \pi\right ) = 0 \Leftrightarrow  \sqrt{-\lambda} =n+\dfrac{1}{2} \, \, (n=0,1,2,\dots) \Leftrightarrow \lambda=-\left(n+\dfrac{1}{2} \right)^2 \, \, (n=0,1,2,\dots),$$
and, in each of these cases, we have the following nontrivial solutions of (\ref{x}):
$$X_n(x)=\cos \left( nx +\dfrac{x}{2} \right) \quad \mbox{and} \quad \tilde{X}_n(x)=\sin \left( nx +\dfrac{x}{2} \right) \, \, (n=0,1,2,\dots).$$
Now we solve (\ref{t}) for $\lambda=-(n+1/2)^2$ and we get the basic nontrivial solutions
$$T_n(t)={\rm e}^{-(n+1/2)^2 k t}, \, \, t \ge 0.$$
Finally, we solve (\ref{edp2})--(\ref{b22}) via superposition principle with the basic solutions $T_n(t)X_n(x)$ and $T_n(t)\tilde{X}_n(x)$ as follows: 
$$v(x,t)=\sum_{n=0}^{\infty} {\rm e}^{-(n+1/2)^2 k t} \left(A_n \cos \left( nx +\dfrac{x}{2} \right)+B_n   \sin \left( nx +\dfrac{x}{2} \right) \ \right),$$
for some adequate constants $A_n, \, B_n \in \mathbb R$ ($n=0,1,2,\dots$) such that (\ref{ic2}) is fulfilled, i.e. for all $x \in [-\pi,\pi]$ we must have 
$$f(x)-1=v(x,0)=\sum_{n=0}^{\infty} \left(A_n  \cos \left( nx +\dfrac{x}{2} \right)+B_n    \sin \left( nx +\dfrac{x}{2} \right) \ \right).$$
We deduce from Theorem \ref{t1} applied to $f(x)-1$ (which is an antiperiodic function) that we must take
$$A_n= \dfrac{1}{\pi}\int_{-\pi}^{\pi}{(f(x)-1) \cos \left(nx+\dfrac{x}{2} \right) \, dx}$$
and
$$B_n= \dfrac{1}{\pi}\int_{-\pi}^{\pi}{(f(x)-1) \cos \left(nx+\dfrac{x}{2} \right) \, dx}.$$
To sum up, according to (\ref{solu}) the solution of (\ref{edp})--(\ref{b2}) is
$$u(x,t)=1+\sum_{n=0}^{\infty} {\rm e}^{-(n+1/2)^2 k t} \left(A_n  \cos \left( nx +\dfrac{x}{2} \right)+B_n    \sin \left( nx +\dfrac{x}{2} \right) \right).$$
In particular, the solution for $k=1$ and $f(x)=(x/\pi)^2$ is
$$u(x,t)=1-\dfrac{32}{\pi^3} \sum_{n=0}^{\infty}\dfrac{(-1)^n}{(2n+1)^3} {\rm e}^{-(n+1/2)^2 t}    \cos \left( nx +\dfrac{x}{2} \right)  ,$$
which tends to $1$ as $t$ tends to infinity; see Figures 7 and 8.
\begin{figure}[h]
\centering
\includegraphics[width=7cm,keepaspectratio=true]{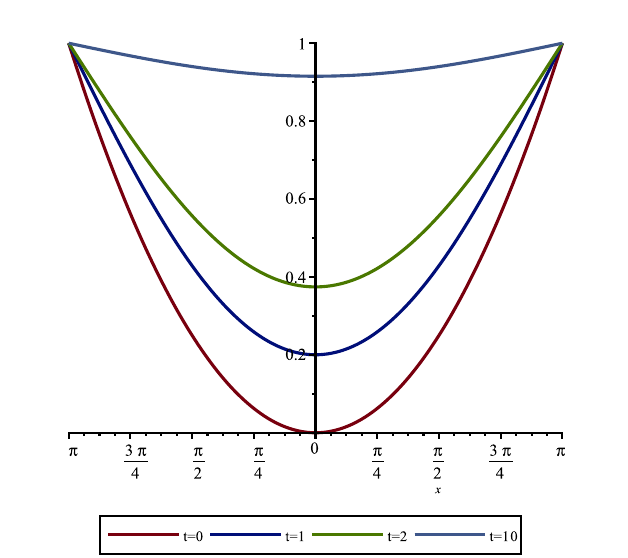}
\caption{Tenth partial sum of $u(x,t)$ for some values of $t$.}
\end{figure}

\begin{figure}[h]
\centering
\includegraphics[width=7cm,keepaspectratio=true]{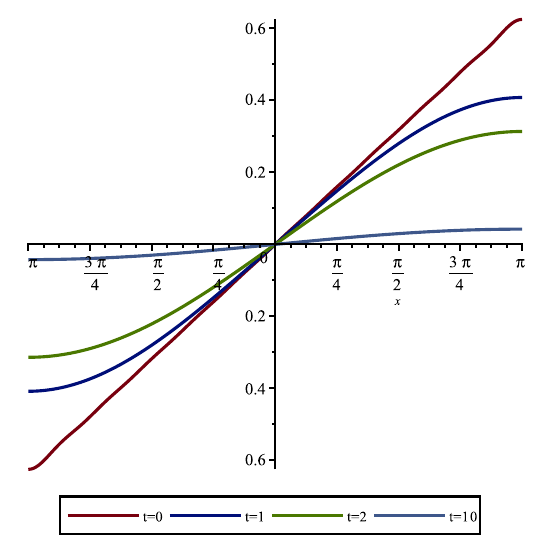}
\caption{Tenth partial sum of $u_x(x,t)$ for some values of $t$.}
\end{figure}

\begin{remark}
All the calculations and graphs in this paper were done with Maple but many other types of software (even open--source software) could have been used instead. Readers are referred to our elementary paper \cite{c}, which contains some useful Sage codes for the Fourier series and their applications.
\end{remark}

\end{document}